\documentclass[11pt, reqno, psamsfonts]{amsart}
\pdfoutput=1

\usepackage{amssymb}
\usepackage{amsthm}
\usepackage{amsmath}
\usepackage{latexsym}
\usepackage[T1]{fontenc}
\usepackage[utf8]{inputenc}
\usepackage[russian, french, english]{babel}

\usepackage{graphicx}
\usepackage{wrapfig}
\usepackage[labelfont=bf]{caption}
\usepackage{mathtools}
\usepackage{tikz}
\usepackage{amsbsy}
\usepackage[inline]{enumitem}
\usepackage{mathrsfs}
\usetikzlibrary{shapes,snakes}
\usetikzlibrary{arrows.meta}
\usetikzlibrary{decorations.pathmorphing}
\usetikzlibrary{patterns}
\usepackage{float}
\usepackage{array}
\usepackage{multicol}
\usepackage{stmaryrd}
\usepackage{eucal,eufrak}

\usepackage{caladea}
\usepackage[T1]{fontenc}
\usepackage[euler-digits,small]{eulervm}
\usepackage[bb=pazo,scr=esstix]{mathalfa}

\usepackage{upgreek}
\usepackage{titlesec}
\usepackage[spacing=true,kerning=true,babel=true,tracking=true]{microtype}
\usepackage[shortcuts]{extdash}
\usepackage[foot]{amsaddr} 
\usepackage{framed}
\usepackage{algpseudocode}
\usepackage{algorithm}
\definecolor{shadecolor}{gray}{0.95}
\usepackage[left=1in,right=1in,top=1in,bottom=1in,bindingoffset=0cm]{geometry}
\usepackage{xspace}
\usepackage{bm}

\usepackage[hidelinks]{hyperref}
\usepackage[backend=biber, style=alphabetic, sorting=nyt, maxnames=100,backref=true]{biblatex}
\title{\lsstyle Local Coloring Problems on Smooth Graphs}
\date{}
\author{\lsstyle Anton~Bernshteyn}
\address{\textls{\normalfont{}School of Mathematics, Georgia Institute of Technology, Atlanta, GA, USA}}
\email{bahtoh@gatech.edu}
\thanks{This research was partially supported by the NSF grant DMS-2045412.}

\newtheoremstyle{bfnote}%
{}{}%
{\slshape}{}%
{\bfseries}{\bfseries.}%
{ }%
{\thmname{#1}\thmnumber{ #2}\thmnote{ \ep{\normalfont{}#3}}}

\newtheoremstyle{defbfnote}%
{}{}%
{}{}%
{\bfseries}{.}%
{ }%
{\thmname{#1}\thmnumber{ #2}\thmnote{ (#3)}}

\newtheoremstyle{claim}%
{}{}%
{\slshape}{}%
{\itshape}{.}%
{ }%
{\thmname{#1}\thmnumber{ #2}\thmnote{ \ep{\normalfont{}#3}}}

\newtheoremstyle{smalldefn}%
{}{}%
{}{}%
{\itshape}{.}%
{ }%
{\thmname{#1}\thmnumber{ #2}\thmnote{ \ep{\normalfont{}#3}}}

\theoremstyle{bfnote}

\newtheorem{theo}[equation]{Theorem}
\newtheorem*{theo*}{Theorem}
\newtheorem{prop}[equation]{Proposition}
\newtheorem{lemma}[equation]{Lemma}

\newtheorem*{claim*}{Claim}
\newtheorem{big_claim}[equation]{Claim}
\newtheorem*{corl*}{Corollary}

\theoremstyle{claim}

\theoremstyle{smalldefn}

\newcommand*{\myproofname}{Proof}

\theoremstyle{definition}
\newtheorem{defn}[equation]{Definition}
\newtheorem*{defn*}{Definition}

\newtheorem*{exmp*}{Example}

\newtheorem*{assum*}{Assumptions}

\theoremstyle{remark}
\newtheorem*{ques*}{Question}
\newtheorem*{remk*}{Remark}

\newcommand{\0}{\varnothing}
\newcommand{\set}[1]{\{#1\}}

\newcommand{\N}{{\mathbb{N}}}
\newcommand{\Z}{\mathbb{Z}}
\newcommand{\F}{\mathbb{F}}
\newcommand{\R}{\mathbb{R}}
\renewcommand{\C}{\mathbb{C}}
\newcommand{\Q}{\mathbb{Q}}

\renewcommand{\epsilon}{\varepsilon}
\renewcommand{\phi}{\varphi}
\renewcommand{\theta}{\vartheta}

\renewcommand{\geq}{\geqslant}

\newcommand{\concat}{{^\smallfrown}}

\renewcommand{\G}{\Gamma}

\newcommand{\defeq}{\coloneqq}

\newcommand{\emphd}[1]{{\fontseries{b}\selectfont{#1}}}

\renewcommand{\P}{\mathbb{P}}

\newcommand{\B}{\mathcal{B}}

\newcommand{\bemph}[1]{{\normalfont#1}} 
\newcommand{\ep}[1]{\bemph{(}#1\bemph{)}} 

\newenvironment{scproof}[1][Proof]{\begin{proof}[#1]}{\end{proof}}

\numberwithin{equation}{section}

\renewcommand{\thesubsection}{\arabic{section}.\Alph{subsection}}

\usepackage{etoolbox}

\titleformat{\section}[block]{\scshape\filcenter}{\thesection.}{1ex}{}
\titleformat{\subsection}[block]{\bfseries\filcenter}{\thesubsection.}{1ex}{}
\titleformat{\subsubsection}[runin]{\itshape}{\bfseries\upshape\thesubsubsection.}{1ex}{}[.---]

\titlespacing*{\section}{0pt}{*3}{*1}
\titlespacing*{\subsection}{0pt}{*3}{*1}
\titlespacing*{\subsubsection}{0pt}{*0.8}{*0}

\makeatletter
\newcommand{\neutralize}[1]{\expandafter\let\csname c@#1\endcsname\count@}
\makeatother

\renewcommand{\mkbibnamefirst}[1]{\textsc{#1}}
\renewcommand{\mkbibnamelast}[1]{\textsc{#1}}
\renewcommand{\mkbibnameprefix}[1]{\textsc{#1}}
\renewcommand{\mkbibnameaffix}[1]{\textsc{#1}}
\renewcommand{\finentrypunct}{}

\renewbibmacro{in:}{}

\renewbibmacro*{volume+number+eid}{%
	\printfield{volume}%
	\setunit*{\addnbspace}
	\printfield{number}%
	\setunit{\addcomma\space}%
	\printfield{eid}}

\DeclareFieldFormat[article]{volume}{\textbf{#1}\space}
\DeclareFieldFormat[article]{number}{\mkbibparens{#1}}

\DeclareFieldFormat{journaltitle}{#1,}
\DeclareFieldFormat[thesis]{title}{\mkbibemph{#1}\addperiod}
\DeclareFieldFormat[article, unpublished, thesis]{title}{\mkbibemph{#1},}
\DeclareFieldFormat[book]{title}{\mkbibemph{#1}\addperiod}
\DeclareFieldFormat[unpublished]{howpublished}{#1, }

\DeclareFieldFormat{pages}{#1}

\DeclareFieldFormat[article]{series}{Ser.~#1\addcomma}

\setlist{topsep=3pt,itemsep=3pt}

\begin{document}
	\pagestyle{plain}
	
	\vspace*{0pt}
	
	\maketitle
	
	\begin{abstract}
		We construct a smooth locally finite Borel graph $G$ and a local coloring problem $\Pi$ such that $G$ has a coloring $V(G) \to \N$ that solves $\Pi$, but no such coloring can be Borel.
	\end{abstract}
	
	\section{Introduction}
	
	By default, graphs in this note are undirected and simple. A graph $G$ is \emphd{Borel} if $V(G)$ is a standard Borel space and $E(G)$ is a Borel subset of $V(G) \times V(G)$. A graph $G$ is \emphd{locally countable} \ep{resp.~\emphd{locally finite}} if every vertex of $G$ has countably \ep{resp.~finitely} many neighbors. Graphs are often allowed to carry additional structure, such as an orientation or a vertex labeling. This idea is captured in the general notion of a ``structured graph.'' The precise definitions are a bit tedious and can be found in \cite[\S2.A.1]{BerDist}. In this note we shall only work with graphs that carry some specific type of structure \ep{e.g., an orientation}, so there is no need to reproduce the general formalism here. Throughout most of this note, the word ``graph'' stands for ``structured graph.'' Of course, when we talk about Borel graphs, the additional structure is understood to be Borel in an appropriate sense.
	
	A locally countable Borel graph $G$ is \emphd{smooth} if it admits a \emphd{Borel transversal}, i.e., a Borel set $T \subseteq V(G)$ of vertices such that every connected component of $G$ contains exactly one vertex from $T$. It is generally understood that smoothness ``trivializes'' Borel combinatorics. For example, if $G$ is a smooth locally finite Borel graph that has a proper $k$-coloring for some $k \in \N$, this automatically implies that $G$ also has a \emph{Borel} proper $k$-coloring \cite[Proposition 1]{ConleyMiller}. In addition to proper $k$-colorings, this is more generally true for solutions to \emph{local coloring problems}. To define this notion, we need some notation first. Let $G$ be a \ep{structured} graph. Given $r \in \N$ and a vertex $x \in V(G)$, we let $[G, x]_r$ denote the isomorphism type of the rooted radius-$r$ ball around $x$ in $G$. A \emphd{coloring} of $G$ is simply a function $f \colon V(G) \to \N$. 
	If $f$ is a coloring of $G$, we call the pair $(G, f)$ a \emphd{colored graph}. Note that colored graphs are a special type of structured graphs. Isomorphisms between colored graphs are required to preserve the coloring. If $(G, f)$ is a colored graph, we write $[G, f, x]_r$ \ep{instead of the more cumbersome $[(G, f), x]_r$} for the isomorphism type of the rooted radius-$r$ ball around $x$ in $(G, f)$.

	
	\begin{defn}
		A \emphd{local coloring problem} is a pair $\Pi = (r, \mathscr{P})$, where $r \in \N$ and $\mathscr{P}$ is a map that sends isomorphism classes of finite labeled graphs to $\set{\mathsf{pass},\mathsf{fail}}$. A \emphd{$\Pi$-coloring} of a locally finite graph $G$ is a function $f \colon V(G) \to \N$ such that $\mathscr{P}([G, f, x]_r) = \mathsf{pass}$ for all $x \in V(G)$.
	\end{defn}

	Informally, a local coloring problem $\Pi = (r, \mathscr{P})$ is a ``rule'' that decides whether a coloring $f$ of $G$ is ``valid'' by looking at the restrictions of $f$ to balls of radius $r$ in $G$. The prototypical example of a local coloring problem is proper coloring, since whether or not a coloring is proper is determined by its restrictions to radius-$1$ balls.
	
	\begin{theo}[{\cite[Theorem 5.23]{Pikh}}]\label{theo:fin}
		Let $G$ be a smooth locally finite Borel graph and let $\Pi$ be a local coloring problem. If $G$ has a $\Pi$-coloring $V(G) \to k$ for $k \in \N$, then $G$ also has a Borel $\Pi$-coloring $V(G) \to k$. 
	\end{theo}

	Theorem~\ref{theo:fin} formally captures the aforementioned intuition that smoothness ``trivializes'' Borel combinatorics. While this statement itself is well-known to experts in the area, it seems that it was first stated and proved in this generality in the survey article \cite{Pikh} by Pikhurko. Below we sketch an alternative proof using the Borel uniformization theorem for Borel sets with compact fibers.
	
	\begin{scproof}[Proof of Theorem~\ref{theo:fin}]
		Let $T$ be a Borel transversal for $G$. For each $x \in T$, let $G_x$ denote the connected component of $x$ in $G$. Note that $f \colon V(G) \to k$ is a $\Pi$-coloring of $G$ if and only if for each $x \in T$, the restriction of $f$ to $V(G_x)$ is a $\Pi$-coloring of $G_x$. By the Feldman--Moore theorem \cite[Theorem 22.2]{AnushDST}, there is a countable sequence of Borel functions $\phi_0$, $\phi_1$, \ldots{} such that for each $x \in T$, $V(G_x) = \set{\phi_n(x) \,:\, n \in \N}$. Let $C_x \subseteq k^\N$ be the set of all functions $f \colon \N \to k$ satisfying the following conditions:
		\begin{enumerate}[label=\ep{\roman*}]
			\item\label{item:1} if $\phi_n(x) = \phi_m(x)$, then $f(n) = f(m)$; and
			\item\label{item:2} the map $V(G_x) \to k \colon \phi_n(x) \mapsto f(n)$ is a $\Pi$-coloring of $G_x$.
		\end{enumerate} 
		Condition \ref{item:1} ensures that the map in \ref{item:2} is well-defined. Since whether or not a coloring is a $\Pi$-coloring is determined by its restrictions to finite subsets, the set $C_x$ is closed in $k^\N$ and hence it is compact. For the same reason, the set $\set{(x, f) \,:\, f \in C_x} \subseteq T \times k^\N$ is Borel. Finally, $C_x \neq \0$ for all $x \in T$ since $G$ has a $\Pi$-coloring $V(G) \to k$. Therefore, by \cite[Theorem 28.8]{KechrisDST}, there is a Borel function $x \mapsto f_x$ that assigns to each $x \in T$ an element $f_x \in C_x$. Now the map $V(G) \to k$ sending each $\phi_n(x)$ to $f_x(n)$ is a desired Borel $\Pi$-coloring of $G$.
	\end{scproof}

	The proof of Theorem~\ref{theo:fin} given above, as well as the one in \cite{Pikh}, crucially relies on the fact that we are looking for colorings using only a finite number of colors. Prompted by a question of Pikhurko \ep{private communication}, we show that this is necessary, as the analog of Theorem~\ref{theo:fin} for colorings with countably many colors fails.
	
	\begin{theo}\label{theo:inf}
		There exist a smooth locally finite Borel graph $G$ and a local coloring problem $\Pi$ such that $G$ has a $\Pi$-coloring $V(G) \to \N$, but no such coloring can be Borel.
	\end{theo}	
	
	The graph $G$ in our construction is not just locally finite, but has maximum degree $3$ \ep{meaning that every vertex of $G$ has at most $3$ neighbors}. It is also acyclic, i.e., all its connected components are trees. 
	
	\section{Proof of Theorem~\ref{theo:inf}}\label{sec:proof}
	
	Recall that a subset $A \subseteq X$ of a standard Borel space $A$ is \emphd{analytic} if it is an image of a Borel set under a Borel function, and \emphd{co-analytic} if $X \setminus A$ is analytic. An analytic \ep{resp.~co-analytic} set $A \subseteq X$ is \emphd{complete analytic} \ep{resp.~\emphd{complete co-analytic}} if for every standard Borel space $Y$ and every analytic \ep{resp.~co-analytic} set $B \subseteq Y$, there is a Borel function $f \colon Y \to X$ such that for all $y \in Y$,
	\[
	y \in B \quad \Longleftrightarrow \quad f(y) \in A.
	\]
	Such a function $f$ is called a \emphd{Borel reduction} from $B$ to $A$.
	
	Natural examples of complete analytic and complete co-analytic sets arise from spaces of trees \ep{see \cite[\S2 and \S27.A]{KechrisDST} for a detailed discussion}. Let $2^{< \infty}$ be the set of all finite sequences of $0$s and $1$s. A \emphd{set-theoretic binary tree} is a nonempty set $T \subseteq 2^{< \infty}$ closed under taking initial segments. We identify a set-theoretic binary tree $T$ with the \ep{graph-theoretic} rooted tree in which the vertices are the sequences in $T$, the empty sequence $\0$ is the \emphd{root}, and the \emphd{parent} of a sequence $a_0 \ldots a_k$ is $a_0 \ldots a_{k-1}$. Note that each vertex $s \in T$ has at most two \emphd{children}, namely $s \concat 0$ and $s \concat 1$ \ep{here $\concat$ denotes concatenation}. When $s \concat 0$ \ep{resp.~$s \concat 1$} is in $T$, we call it the \emphd{left} \ep{resp.~\emphd{right}} child of $s$. A tree $T$ is \emphd{pruned} if every $s \in T$ has at least one child. For a pruned binary tree $T$, let $[T] \subseteq 2^\N$ be the set of all \emphd{branches} in $T$, i.e., all infinite binary sequences $\beta \in 2^\N$ such that every finite initial segment of $\beta$ is in $T$. The set $\mathrm{PTr}_2$ of all pruned set-theoretic binary trees is a closed subset of $2^{2^{< \infty}}$, and as such it is a standard Borel space. Let $F \subseteq \mathrm{PTr}_2$ denote the set of all pruned binary trees $T$ such that every branch $\beta \in [T]$ contains finitely many $1$s \ep{i.e., $\beta_n =0$ for all but finitely many $n \in \N$}.
	
	\begin{prop}[{\cite[Exercise~27.3]{KechrisDST}}]\label{prop:complete}
		The set $F \subseteq \mathrm{PTr}_2$ is complete co-analytic.
	\end{prop}

	The combinatorial core of our argument is the observation that there exists a local coloring problem $\Sigma$ with the following property: a pruned binary tree $T$ admits a $\Sigma$-coloring if and only if $T \not\in F$. Namely, we define $\Sigma$ so that $f \colon T \to \N$ is a $\Sigma$-coloring if and only if it satisfies the following constraints:
	\begin{enumerate}[label=\ep{$\Sigma$\arabic*}]
		\item\label{item:root} If $x$ is the root of $T$, then $f(x) \geq 1$.
		
		\item\label{item:right} If $x$ is a vertex such that $f(x) = 1$, then $x$ must have a right child $y$ and $f(y) \geq 1$.
		
		\item\label{item:left} If $x$ is a vertex such that $f(x) \geq 2$, then $x$ must have a left child $y$ with $f(y) = f(x) - 1$.
	\end{enumerate}
	Note that this is indeed a local coloring problem, as the above conditions are only concerned with the relationship between the colors of each vertex and its neighbors. 

	\begin{lemma}\label{lemma:red}
		A pruned binary tree $T$ admits a $\Sigma$-coloring if and only if $T \not\in F$.
	\end{lemma}
	\begin{scproof}
		Suppose first that $T \not\in F$. This means that there is a branch $\beta \in [T]$ that contains infinitely many $1$s. For a finite initial segment $s$ of $\beta$, let $k(s)$ denote the \ep{unique} $k \in \N$ such that
		\[
		s \concat \underbrace{0 \ldots 0}_{\text{$k$ zeroes}} \concat 1
		\]
		is an initial segment of $\beta$. Now define $f$ by setting $f(s) \defeq k(s) + 1$ for every finite initial segment $s$ of $\beta$ and sending every other vertex of $T$ to $0$. 
		It is easy to see that $f$ is a $\Sigma$-coloring.
		
		Conversely, suppose $T$ has a $\Sigma$-coloring $f$. Call $s \in T$ \emphd{positive} if $f(s) \geq 1$. For each positive $s \in T$, let the \emphd{favorite child} $\mathrm{fav}(s)$ of $s$ be its right child if $f(s) = 1$ and left child if $f(s) \geq 2$. Conditions \ref{item:right} and \ref{item:left} ensure that every positive $s \in T$ has a favorite child and that $\mathrm{fav}(s)$ is also positive. By \ref{item:root}, the root of $T$ is positive, so we can build a branch $\beta \in [T]$ by iteratively applying the operation $\mathrm{fav}$ to the root of $T$. Note that if $\mathrm{fav}(s)$ is the left child of $s$, then, by \ref{item:left}, $f(\mathrm{fav}(s)) = f(s) - 1$. Since a strictly decreasing sequence of natural numbers cannot be infinite, $\beta$ must contain infinitely many $1$s, and hence $T \not\in F$. 
	\end{scproof}


	For our construction, we shall use the following fact from descriptive set theory: 
	
	\begin{prop}[{\cite[\S35.A]{KechrisDST}}]\label{prop:insep}
		Let $X$ be an uncountable standard Borel space. Then there exist disjoint co-analytic sets $A_0$, $A_1 \subseteq X$ that are \emphd{Borel-inseparable}, meaning that it is impossible to partition $X$ into two Borel sets one of which contains $A_0$ and the other one $A_1$.
	\end{prop}

	Let $X$ be an uncountable standard Borel space and let $A_0$, $A_1 \subseteq X$ be two Borel-inseparable co-analytic subsets of $X$, as given by Proposition~\ref{prop:insep}. Since the set $F$ is complete co-analytic by Proposition~\ref{prop:complete}, there are Borel maps $x \mapsto T_x^0$, $x \mapsto T_x^1$ sending points $x \in X$ to pruned binary trees such that for each $i \in \set{0,1}$,
	\begin{equation}\label{eq:red}
		x \in A_i \quad \Longleftrightarrow \quad T_x^i \in F.
	\end{equation}
	For each $x \in X$ and $i \in \set{0,1}$, define $\hat{T}_x^i \defeq \set{(x, i, s) \,:\, s \in T_x^i}$ and then turn $\hat{T}^i_x$ into a rooted binary tree in the obvious way, so that $T_x^i \cong \hat{T}_x^i$ via the mapping $s \mapsto (x, i, s)$. In particular, the root of $\hat{T}^i_x$ is $(x, i,\0)$. By construction, the trees $\hat{T}_x^i$ are pairwise disjoint, and their union is a Borel graph $H$ with vertex set
	\[
		V(H) \,\defeq\, \set{(x, i, s) \,:\,x \in X, \  i \in \set{0,1}, \ \text{and}\ s \in T_x^i}.
	\]
	Without loss of generality, we may assume that $X \cap V(H) = \0$ and form a graph $G$ with vertex set $V(G) \defeq V(H) \sqcup X$ by adding to $H$ the edges between $x$ and $(x, i, \0)$ for all $x \in X$ and $i \in \set{0,1}$. The vertices $x \in X$ are called the \emphd{anchor vertices} of $G$. 
	Each connected component of $G$ contains exactly one anchor vertex $x$ and looks like a disjoint union of $T_x^0$ and $T_x^1$ with roots connected to $x$ by edges. Observe that $G$ is locally finite \ep{in fact, every vertex of $G$ has at most $3$ neighbors} and smooth \ep{since $X$ is a Borel transversal for $G$}.
	
	Finally, we let $\Pi$ be a local coloring problem such that $f \colon V(G) \to \N$ is a $\Pi$-coloring if and only if:
	\begin{enumerate}[label=\ep{$\Pi$\arabic*}]
		\item\label{item:anch} Every anchor vertex $x$ has exactly one neighbor $y$ with $f(y) \geq 1$.
		
		\item\label{item:right1} If $x$ is a non-anchor vertex such that $f(x) = 1$, then $x$ must have a right child $y$ and $f(y) \geq 1$.
		
		\item\label{item:left1} If $x$ is a non-anchor vertex such that $f(x) \geq 2$, then $x$ must have a left child $y$ with $f(y) = f(x) - 1$.
	\end{enumerate}
	As in the definition of $\Sigma$, these conditions are only concerned with the relationship between the colors of each vertex and its neighbors, so this is indeed a local coloring problem. We remark that $\Pi$ treats $G$ as a {structured} graph that retains the rooted tree structure on each $\hat{T}_x^i$ and also includes a label for each vertex indicating whether it is an anchor vertex. However, it is not hard to modify the construction to ``encode'' all the necessary extra information in the graph structure; we explain how to do this in \S\ref{sec:no_structure}.
	
	
	\begin{big_claim}
		The graph $G$ has a $\Pi$-coloring.
	\end{big_claim}
	\begin{scproof}
		It is enough to argue that every connected component of $G$ has a $\Pi$-coloring. To this end, consider any anchor vertex $x \in X$. Since $A_0$ and $A_1$ are disjoint, there is $i \in \set{0,1}$ such that $x \not \in A_i$. By \eqref{eq:red}, this means that $T_x^i \not \in F$, and hence, by Lemma~\ref{lemma:red}, $T_x^i$ has a $\Sigma$-coloring $h \colon T_x^i \to \N$. Thus, we can define a $\Pi$-coloring $f$ on the connected component of $x$ by setting $f(x, i, s) \defeq h(s)$ for all $s \in T_x^i$ and sending every other vertex to $0$.
	\end{scproof}

	\begin{big_claim}
		The graph $G$ has no Borel $\Pi$-coloring.
	\end{big_claim}
	\begin{scproof}
		Suppose, toward a contradiction, that $f \colon V(G) \to \N$ is a Borel $\Pi$-coloring. For each $i \in \set{0,1}$, let
		\[
			B_i \,\defeq\, \set{x \in X \,:\, f(x, i, \0) \geq 1}.
		\]
		By \ref{item:anch}, $X = B_0 \sqcup B_1$ is a partition of $X$, and since $f$ is a Borel function, the sets $B_0$ and $B_1$ are Borel. We claim that $B_0 \cap A_0 = B_1 \cap A_1 = \0$. Indeed, take any $x \in B_i$. Since $f(x,i,\0) \geq 1$, conditions \ref{item:right1} and \ref{item:left1} imply that the restriction of $f$ to $\hat{T}_x^i$ is a $\Sigma$-coloring. By Lemma~\ref{lemma:red}, this means that $T_x^i \not \in F$, and hence $x \not \in A_i$ by \eqref{eq:red}, as claimed. Since $X = B_0 \sqcup B_1$, 
		we conclude that $A_0 \subseteq B_1$ and $A_1 \subseteq B_0$, which 
		contradicts the fact that $A_0$ and $A_1$ are Borel-inseparable.
		%
	\end{scproof}

	\section{Removing the extra structure}\label{sec:no_structure}
	
	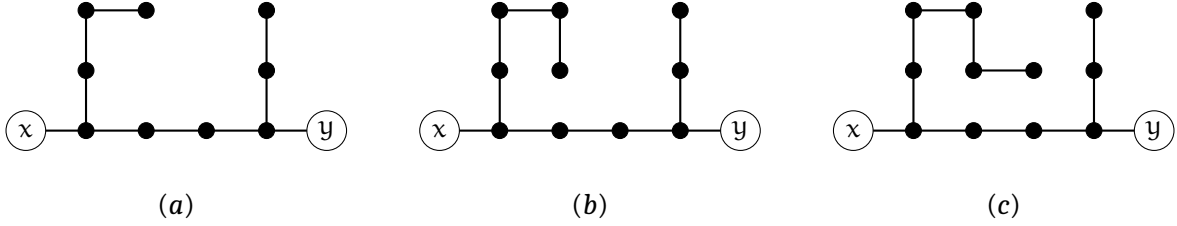
\begin{figure}[t]
		\centering
		\begin{tikzpicture}
		\begin{scope}
			\node[circle, draw, minimum size=15pt,outer sep=0pt,inner sep=0pt] (x) at (0,0) {$x$};
			\node[circle, draw, minimum size=15pt,outer sep=0pt,inner sep=0pt] (y) at (4,0) {$y$};	
			\filldraw (0.8,0) circle (3pt);
			\filldraw (1.6,0) circle (3pt);	
			\filldraw (2.4,0) circle (3pt);	
			\filldraw (3.2,0) circle (3pt);	
			
			\filldraw (3.2,0.8) circle (3pt);	
			\filldraw (3.2,1.6) circle (3pt);
			\draw[thick] (3.2,0) -- (3.2,1.6);
			
			\filldraw (0.8,0.8) circle (3pt);	
			\filldraw (0.8,1.6) circle (3pt);
			\filldraw (1.6,1.6) circle (3pt);
			\draw[thick] (0.8,0) -- (0.8,1.6) -- (1.6,1.6);
			
			\draw[thick] (x) -- (y);	
			\node at (2,-1) {\ep{\emph{a}}};						
		\end{scope}
		
		\begin{scope}[xshift=5.5cm]
		\node[circle, draw, minimum size=15pt,outer sep=0pt,inner sep=0pt] (x) at (0,0) {$x$};
		\node[circle, draw, minimum size=15pt,outer sep=0pt,inner sep=0pt] (y) at (4,0) {$y$};	
		\filldraw (0.8,0) circle (3pt);
		\filldraw (1.6,0) circle (3pt);	
		\filldraw (2.4,0) circle (3pt);	
		\filldraw (3.2,0) circle (3pt);	
		
		\filldraw (3.2,0.8) circle (3pt);	
		\filldraw (3.2,1.6) circle (3pt);
		\draw[thick] (3.2,0) -- (3.2,1.6);
		
		\filldraw (0.8,0.8) circle (3pt);	
		\filldraw (0.8,1.6) circle (3pt);
		\filldraw (1.6,1.6) circle (3pt);
		\filldraw (1.6,0.8) circle (3pt);
		\draw[thick] (0.8,0) -- (0.8,1.6) -- (1.6,1.6) -- (1.6,0.8);
		
		\draw[thick] (x) -- (y);	
		\node at (2,-1) {\ep{\emph{b}}};									
		\end{scope}
		
		\begin{scope}[xshift=11cm]
		\node[circle, draw, minimum size=15pt,outer sep=0pt,inner sep=0pt] (x) at (0,0) {$x$};
		\node[circle, draw, minimum size=15pt,outer sep=0pt,inner sep=0pt] (y) at (4,0) {$y$};	
		\filldraw (0.8,0) circle (3pt);
		\filldraw (1.6,0) circle (3pt);	
		\filldraw (2.4,0) circle (3pt);	
		\filldraw (3.2,0) circle (3pt);	
		
		\filldraw (3.2,0.8) circle (3pt);	
		\filldraw (3.2,1.6) circle (3pt);
		\draw[thick] (3.2,0) -- (3.2,1.6);
		
		\filldraw (0.8,0.8) circle (3pt);	
		\filldraw (0.8,1.6) circle (3pt);
		\filldraw (1.6,0.8) circle (3pt);
		\filldraw (1.6,1.6) circle (3pt);
		\filldraw (2.4,0.8) circle (3pt);
		\draw[thick] (0.8,0) -- (0.8,1.6) -- (1.6,1.6) -- (1.6,0.8) -- (2.4,0.8);
		
		\draw[thick] (x) -- (y);	
		\node at (2,-1) {\ep{\emph{c}}};								
		\end{scope}
		\end{tikzpicture}
		\caption{The ``gadgets'' replacing the edges of $G$ in $G^\ast$. If $x$ is an anchor vertex and $y$ is its neighbor, then the edge $xy$ is replaced by the graph shown in \ep{\emph{a}}. If $x$ and $y$ are non-anchor vertices and $y$ is a left \ep{resp. right} child of $x$, then the edge $xy$ is replaced by the graph shown in \ep{\emph{b}} \ep{resp.~\ep{\emph{c}}}.}\label{fig:gadgets}
	\end{figure}
	
	In this section we sketch how the construction from \S\ref{sec:proof} can be modified to make $G$ a ``plain'' graph without any extra structure. Recall that, in addition to $G$'s adjacency relation, the problem $\Pi$ uses the following information:
	\begin{itemize}
		\item which vertices of $G$ are the anchor vertices;
		
		\item which one of each pair of adjacent non-anchor vertices is the parent and which one is the child;
		
		\item if $y$ is a child of $x$, whether $y$ is the left or the right child.
	\end{itemize}
	All this can be encoded by replacing each edge of $G$ with one of the three finite ``gadgets'' shown in Fig.~\ref{fig:gadgets}. 
	The resulting graph $G^\ast$ still is acyclic and has maximum degree $3$. We view $V(G)$ as a subset of $V(G^\ast)$; in other words, the vertices of $G^\ast$ are of two types: the \emphd{original} vertices of $G$ and the new \emphd{auxiliary} vertices \ep{the auxiliary vertices are shown in Fig.~\ref{fig:gadgets} in black}. Note that every Borel transversal for $G$ is also a Borel transversal for $G^\ast$, so $G^\ast$ is smooth. 
	Say that a vertex $x \in V(G^\ast)$ has \emphd{order} $k$ if $G^\ast$ contains a path $x$---$x_1$---$\cdots$---$x_k$ such that the vertices $x_1$, \ldots, $x_{k-1}$ have degree $2$ and the vertex $x_k$ has degree $1$. Let us now explain how to recover the structure of $G$ from $G^\ast$.
	\begin{itemize}
		\item A vertex $x \in V(G^\ast)$ is original if and only if all its neighbors have degree $3$. 
		
		\item An original vertex is an anchor vertex if and only if it has a neighbor of order $3$.
		
		\item Two original vertices are adjacent in $G$ if and only if they are joined by a path of length $5$ in $G^\ast$.
		
		\item Let $x$ and $y$ be adjacent non-anchor vertices in $G$ and let the $xy$-path in $G^\ast$ be $x$---$u$---$v$---$w$---$z$---$y$. Then $y$ is a left \ep{resp.~right} child of $x$ if and only if $u$ has order $4$ \ep{resp.~$5$}.
	\end{itemize}
	Hence, by only using the graph structure of $G^\ast$, we can define a local coloring problem $\Pi^\ast$ such that $f \colon V(G^\ast) \to \N$ is a $\Pi^\ast$-coloring of $G^\ast$ if and only if the restriction of $f$ to $V(G)$ is a $\Pi$-coloring of $G$. Then the graph $G^\ast$ has a $\Pi^\ast$-coloring but not a Borel one, as desired.
	
	\subsubsection*{{Acknowledgments}}
	
	I am grateful to Oleg Pikhurko for asking the question that prompted this work and for insightful discussions, to Anush Tserunyan for useful comments, and to the anonymous referee for helpful suggestions.
	
	\printbibliography

\end{document}